\input AHTOHFIE.STY
\UDC{\hskip-0.1em
519.157.1
+519.175.1
+519.176
}

\def\F{{\cal F}}
\def\repV{\Upsilon_{\rm v}}
\def\repVSym{\Upsilon_{\rm v}^{\rm sym}}
\def\repE{\Upsilon_{\rm e}}
\def\repESym{\Upsilon_{\rm e}^{\rm sym}}

\title{
Invariant systems of representatives,
\\
or
\\
The cost of symmetry
}

\author{%
Anton A. Klyachko$^\sharp$
\quad
Natalia M. Luneva$^\flat$
}
\address{
$^\sharp$
Faculty of Mechanics and Mathematics of Moscow State University,
Moscow 119991, Leninskie gory, MSU.\\
\strut
$^\flat$%
Sberbank of Russia,
Moscow 117997, ul. Vavilova, 19.
\\
klyachko@mech.math.msu.su
\qquad\qquad\qquad
ssl.nattie@gmail.com
}

\grantsFirst{\RFBR19-01-00591}

\abstract{%
\narrower
Suppose that one can destroy all 
$k$-gons
in a graph
by removing
$n$ 
edges.
How many edges must be removed to destroy all
$k$-gons 
in such a way that the set of removed edges is
invariant with respect to all automorphisms the initial graph?
This paper contains solutions to such kind of problems.
Several open questions are raised.}

\s 0.
Introduction

Consider the following ``applied'' problem.

\disp{\narrower%
We are recruiting a team for the Mars expedition, and we want to
satisfy (for example) the following \emph{compatibility requirement:}
{among any five participants, there should be at least two,
each of
which respects at least three members of
this five.}
Our dossiers show that the expulsion of ten
particular candidates would make this
requirement satisfied.
The problem is that we want to be
\emph{fair {\rm and} impartial}, i.e.,
{we want the set of expelled candidates
to be invariant under all
permutations
\(of all candidates\)
preserving the relation ``respects''.}\
How many
candidates do we need to expel (in the worst case)?
}%
What is the cost of fairness?
The question is not very trivial. For instance, if we try to expel
all candidates obtained from the initial ten ``bad''
candidates by the action of all permutations preserving the relation
``respects", then we can end up with expelling \emph{all}
candidates, even if there are infinitely many of them. In fact, the
optimal set of fairly expelled candidates is always finite and needs
neither contain the given set of ten ``bad'' candidates nor be
contained in it.

The mathematical form of this problem is the following:
\disp{%
Given a family $\F$ of directed five-vertex graphs 
(\emph{forbidden graphs}) and a directed graph~$\Gamma$.
Let $c(\F,\Gamma)$ be the minimal cardinality of
a set of ``marked" vertices such that 
\item{1)}
each subgraph of $\Gamma$ isomorphic to a graph from $\cal F$
contains a marked vertex;
\enditem
and let $c_{\rm sym}(\F,\Gamma)$
be the minimal cardinality of
a set of ``marked" vertices 
satisfying Condition 1) and the following Condition 2): 
\item{2)}
the set of marked vertices is invariant with respect to 
all automorphisms of $\Gamma$.
\enditem
How can $c_{\rm sym}(\F,\Gamma)$ can be bounded in terms of
$c(\F,\Gamma)$\?
}%
\noindent
For this particular Mars-expedition problem,
$$
\eqalign{
\F=&\{\hbox{5-vertex digraphs 
containing at most one vertex
of outdegree $\ge3$}\};
\cr
&\hbox{\quad
how to estimate 
$c_{\rm sym}(\F,\Gamma)$
given that
$c(\F,\Gamma)\le10$\?}
}
$$

\noindent
In algebra, there are many theorems of this kind, e.g.,
\-
if a group $G$ contains an abelian subgroup of finite index, then
$G$ contains a \emph{characteristic} (i.e. invariant with respect to
all automorphisms)
abelian subgroup of finite index [KaM82];
\-
if a group $G$ contains a nilpotent subgroup of finite index, then
$G$ contains a characteristic
nilpotent
(of the same class)
subgroup of finite index [BrNa04];
\-
if a group $G$ contains a solvable subgroup of finite index, then
$G$ contains a characteristic
solvable (of the same derived length) subgroup of finite
index [KhM07a];
\-
if a group $G$ contains a central metabelian
subgroup
of finite index, then
$G$ contains a characteristic
central metabelian
subgroup of finite index [KhM07a];
\-
if a group $G$ contains a paranilpotent subgroup of finite index, then
$G$ contains a characteristic
paranilpotent subgroup of finite index [dGT19b];
\-
if a group $G$ contains a finite-index subgroup
whose commutator subgroup is finite, then
$G$ contains a characteristic
finite-index subgroup whose commutator subgroup is finite [KlMi15];
\-
if a group $G$ of finite exponent contains a finite
normal subgroup $N$, then $G$ contains a characteristic
finite subgroup $H$ such that the \emph{spectrum} (i.e. set of orders of
all elements) of the
quotient group $G/H$ is contained in the spectrum of $G/N$
[KlMi15];
\-
if an algebra $G$ (associative or Lie)
over a field
contains a solvable
ideal of finite codimension, then $G$ contains
an invariant with respect to all automorphisms
solvable
(of the same derived length) ideal of finite codimension~[KhM08].

\enditem
The list can be extended, see, e.g.,
[Vd00],
[KhM07b],
[KhM08],
[KlMe09],
[KhKMM09],
[MSh12],
[KlMi15],
[Fr18],
[dGT18a],
[dGT18b],
[dGT19a],
[dGT19b],
and references therein.
(There are probably similar facts from other branches of mathematics, 
but unfortunately we cannot remember any non-algebraical ones.)  
In algebra, the
assertions of this kind are called sometimes
([KlMi15], [Fr18])
\emph{Khukhro--Makarenko type theorems}
after one of such results [KhM07a] (see also [KlMe09])
including as special cases the first four of the
listed facts.

The first of these facts
(on abelian subgroups) is simple in the sense that there is a short and
elementary proof. However, quite na\"\i ve approaches do not
work. How can one make a given abelian finite-index subgroup~$N\subset G$
characteristic?
\-
One can take the intersection of all automorphic images of~$N$.
The obtained subgroup $\bigcap\limits_{\phi\in\Aut G}\phi(N)$ is
characteristic and abelian (because it is contained in $N$),
but, unfortunately, the index can be infinite.
\-
One can take the subgroup generated by all automorphic
images of $N$. The obtained subgroup
$\gp{\bigcup\limits_{\phi\in\Aut G}\phi(N)}$ is characteristic and
of finite index (because it contains
$N$), but, unfortunately, it can be non-abelian.

\enditem
Actually,
it may happen that
the sought characteristic abelian finite-index subgroup
neither contain nor is contained in $N$.
The same is true for
other Khukhro--Makarenko type theorems (and for the Mars-expedition
problem too).

Almost all the listed
Khukhro--Makarenko type theorems
and the Khukhro--Makarenko theorem itself
are special cases of a
very
general fact [KlMi15]
referred hereafter as the
\emph{multilinear-property theorem}.
This general assertion
gives also a universal
estimate of the corresponding
parameter (the index of the characteristic subgroup, the codimension of
the automorphism-invariant ideal and so on). This
estimate can however be far from sharp in particular cases. For
instance, the multilinear-property theorem gives
the following quantitative refinement of the first listed fact:
\disp{\sl\hfuzz15pt
if a group $G$
contains an abelian subgroup of a finite index $n$, then $G$ contains
a characteristic abelian subgroup of index at most
$(n!)^{\log_2(n!)+1}$.
}%
While the sharp estimate here is $n^2$, see [PSz02].%
\fn{%
See also paper
[dGT18a], whose
authors reproved the estimate $n^2$ apparently
being unaware of~[PSz02]. For finite groups $G$, this
estimate was known earlier ([ChD89], see also [Is08], Theorem~1.41).
While the qualitative fact
-- {\sl if a group contains an abelian finite-index subgroup, then
it contains a {characteristic}
abelian finite-index subgroup} -- was known even earlier,
see~[KaM82].
}
In principle, the
multilinear-property theorem is applicable also to purely combinatorial
questions (see [KlMi15]);
e.g., for the Mars-expedition problem, this theorem gives the following
practically useless answer:
\disp{\sl
it is sufficient to expel
at most
$22 229 709 804 712 410=f(f(f(f(10))))$
candidates, where $f(x)=x(x+1)$.
}%
This bound (greatly exceeding the Earth's population)
is also very coarse. The following theorem says that
the actual answer is 50 (and this estimate is sharp).

\proclaim{Main theorem}.
Suppose that a group $G$ acts on a set $U$, and
$\F$ is a $G$-invariant family
of finite
subsets of
$U$
of
uniformly bounded
cardinality \(i.e. $\max\limits_{F\in\F}|F|<\infty$\).
Let $X\subseteq U$ be a
finite
system of representatives for this
family \(i.e. $X\cap F\ne\emptyset$ for any~$F\in\F$\).
Then there exists a $G$-invariant system of representatives $Y$ such that
$
|Y|\le |X|\cdot
\max
\limits_{F\in\F}|F|.
$

Here, the word \emph{family} means an \emph{unordered
family}, i.e., $\F$ is a set of subsets of $U$.
The \emph{invariance} of $\F$
means that
$g\F=\F$ for all $g\in G$, i.e.,
$gF\:=\{gf\;|\; f\in F\}\in\F$ for all $g\in G$ and $F\in\F$.

The proof of the main theorem is elementary, except
that we use a theorem of B. Neumann [Neu54]
on covering groups by cosets (see the last
section).
The main theorem immediately implies the following fact about graphs.

\Corollary 1.
Let $\Gamma$ be a graph and let $K$ be a finite graph. Then
\item{\rm1)}
if $\Gamma$ contains a finite set of vertices $X$ such that each
subgraph of $\Gamma$ isomorphic to $K$ has at least one vertex
from $X$, then $\Gamma$ contains a finite set of vertices $Y$
invariant with respect to all automorphisms of~$\Gamma$ and such that
again each subgraph of~$\Gamma$ isomorphic to $K$ has at
least one vertex from $Y$\; moreover,
$
|Y|\le |X|\cdot(\hbox{\rm the number of vertices of $K$});
$
\item{\rm2)}
if $\Gamma$ contains a finite set of edges $X$ such that each
subgraph of $\Gamma$ isomorphic to $K$ has at least one edge
from $X$, then $\Gamma$ contains a finite set of edges $Y$
invariant with respect to all automorphisms of~$\Gamma$ and such that
again each subgraph of~$\Gamma$ isomorphic to $K$ has at
least one edge from $Y$\; moreover,
$
|Y|\le |X|\cdot(\hbox{\rm the number of edges of $K$});
$

\noindent
The word \emph{graph} hereafter
can be understood in any reasonable sense:
\-
a graph can be directed, undirected, or mixed,
\-
multiple edges and/or loops can be allowed or not allowed.
\enditem
Of course,
the word ``isomorphism''
(and ``automorphism'') should be understood correspondingly, i.e.
an isomorphism must preserve the direction of edges for directed graphs.

\Remark.
An analogue of Corollary 1 holds when there are
finitely many ``forbidden'' finite graphs
$K_1,\dots,K_n$ (instead the single graph $K$). In this case, the
inequality is
$
|Y|\le |X|\cdot\max\limits_i(\hbox{the number of vertices [edges] of
$K_i$}).
$

The proof is easy: to obtain the assertion (about
vertices)
we apply the main theorem
putting
$$
U=\{\hbox{vertices of } \Gamma\},
\quad
G=\Aut(\Gamma),
\quad
\F=\Bigl\{
\{\hbox{vertices of $S$}\}\;\bigm|\;
\hbox{$S$ is a subgraph of $\Gamma$ isomorphic to one of $K_i$}
\Bigr\}.
$$

In Sections 1 and 2, we discuss the sharpness of estimates
from Corollary 1; the situation is the following:
\-
the both estimates
from Corollary 1 are sharp in the sense
that we cannot
replace the functions
\newline
$|X|\cdot(\hbox{number vertices [edges] graph $K$})$ with smaller
functions of $|X|$ and the number vertices [edges] of $K$;
\-
if we fix a graph $K$ and ask about the sharpness of the estimates
for this given $K$, then the problem becomes more interesting:
\itemitem{-}
the class of graphs $K$ such that
the estimate from
Corollary 1(1)
(about vertices) is sharp have a simple description
(but some interesting questions remain open nevertheless),
see the next section;
\itemitem{-}
as for Corollary 1(2)
(about edges),
we have more questions than
answers; see Section~2.

\enditem
In conclusion, note that
there are situations where no Khukhro--Makarenko type theorems
can be proved. For example,
\disp{\sl
the existence of a finite-index subgroup \(in a group\) with
the law $x^{2019}=1$ does not imply the existence
of a characteristic
finite-index subgroup satisfying this law \rm[KhKMM09].
}%
Paper [KlMi15] contains an
amusing example of a boundary
situation:
on the one hand,
\disp{\sl
if a graph can be made planar by
removing finitely many edges, then such a finite set of edges
can be chosen invariant with respect to all automorphisms
of the graph\;
}%
on the other hand,
\emph{no estimate is possible} (i.e.
there exists $n$ such that, for any $m$,
there exists a graph that can be made planar by removing $n$ edges,
but cannot be made planar by removing a
set of edges consisting of less than $m$ edges and
invariant with respect to all automorphisms of the graph).

{\noindent \bf Our Notation and conventions}
are mainly standard. Note only that
the word \emph{graph} have,
except where otherwise indicated,
six meanings specified above.
The word \emph{subgraph} means an arbitrary subgraph (not necessarily 
induced). The \emph{embedding} of graphs is understood similarly.  
%
The index of a subgroup $H$ of a group
$G$ is denoted as $|G{:}H|$. The letter~$\Z$ denotes the set of integers.
The symbol $|X|$ means the cardinality of a set $X$.

\s 1.
Vertex representativeness

We say that the \emph{vertex representativeness} $\repV(K,\Gamma)$
of a graph $K$ in a graph $\Gamma$
is the minimal integer $n$ such that
$\Gamma$ contains a set $X$
of $n$ vertices
satisfying the following condition:
$$
\hbox{\sl each subgraph of $\Gamma$ isomorphic to $K$ contains a vertex
from $X$.}
\eqno{(*)}
$$
Let us define the
\emph{symmetric vertex representativeness}
$\repVSym(K,\Gamma)$ of a graph $K$ in a graph $\Gamma$ as the minimal
integer $n$ such that $\Gamma$ contains an $(\Aut\Gamma)$-invariant
set $X$ of $n$ vertices
satisfying $(*)$.

For instance,
for the faces of the tetrahedron
and cube, we have:

\bigskip

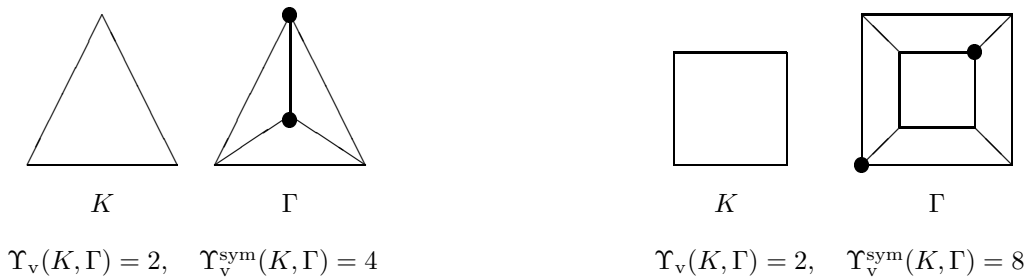
\nobreak
\centerline{\input 0.PIC}
\goodbreak

\bigskip

Clearly, $\repV(K,\Gamma)\le\repVSym(K,\Gamma)$.
Corollary 1 says that
$
\repVSym(K,\Gamma)\le\repV(K,\Gamma)\cdot(\hbox{the number of vertices
of $K$})
$.
The following assertion shows that, for a connected graphs~$K$,
this estimate is sharp.

A graph $K$ is called
\emph{costly in the sense of symmetric vertex representativeness}
or simply
\emph{costly} (or \emph{vertex-costly})~if
$$
\hbox{
$\forall m\in\Z$ there exists a graph
$\Gamma_m$ such that
}
\repVSym(K,\Gamma_m)=
\repV(K,\Gamma_m)\cdot(\hbox{the number of vertices graph $K$})\ge m.
\eqno{(**)}
$$
Thus, $K$ is costly if the estimate from corollary 1(1)
is sharp
for~$K$.
We call a costly graph $K$ \emph{\(vertex-\)costly in a class of
graphs $\cal K$} if the graphs $\Gamma_m$ in $(**)$
can chosen from the class~$\cal K$.

\Th 1.
A finite graph $K$ is costly
if and only if it is connected.
Moreover, any connected graph $K$ without hanging edges
is costly in the class of
connected graphs.

\Proof
Note that
\disp{\sl
any graph $K$ embeds into a vertex-transitive graph
$\~K$ with the same number of vertices.
}%
Indeed, if the word \emph{graph} means an undirected
graph without multiple edges and loops, then we can take the complete
graph as $\~K$; in other cases,
this fact remains valid (we leave it to readers as
an exercise, see graphs $K$ and $\~K$ in Figure~1).

If $K$ is connected, then
we can take the disjoint
union of $m$ copies of~$\~K$ as the graph $\Gamma_m$.
Indeed, to represent all subgraphs isomorphic to $K$
it is
sufficient (and necessary)
to take one vertex from each copy of
$\~K$. Thus,
$\repV(K,\Gamma_m)=m$.
The graph $\Gamma_m$ is vertex-transitive, therefore,
$\repVSym(K,\Gamma_m)=mk$, where $k$ is the number of vertices of $K$.

To prove the assertion ``Moreover'', we add to
this graph $\Gamma_m$ chains of length~$N$ joining each
vertex of the $i$-th copy of~$\~K$ with the corresponding vertex
of the $(i+1)$-th
copy of~$\~K$, where $i\in\{1,\dots,m-1\}$ and the integer~$N$ (the
same for all chains) is larger than the number of vertices of~$K$,
see Figure 1.


\goodbreak
\bigskip
\centerline{\input 1.PIC}
\nobreak%
\centerline{Fig. \lowercase{1}}%
\goodbreak
\bigskip

Again, $\repV(K,\Gamma_m)=m$, because to represent all subgraphs
isomorphic to~$K$, we should mark one vertex
from each copy of~$\~K$ (because the graph $K$ has no hanging
edges and the number $N$ is sufficiently large);
$\repVSym(K,\Gamma_m)=mk$, because
the group $\Aut\Gamma_m$ acts transitively on
vertices of each copy of~$\~K$.

It remains to show that no disconnected graph $K$ is costly.
Let us prove slightly more:
$$
\repVSym(K,\Gamma)\le k_1(\repV(K,\Gamma)+k_2)
$$
if $K=K_1\sqcup K_2$
where  $k_1$ is the number of vertices of $K_1$, and $k_2\le k_1$ is the
number vertices of $K_2$.

Let us mark $\repV(K,\Gamma)$ vertices in graph $\Gamma$ in such a way
that each subgraph isomorphic to $K=K_1\sqcup K_2$ have a marked
vertex. If $\Gamma$ has a subgraph~$\^K_2\iso K_2$ without
marked vertices, then it must intersect each subgraph
isomorphic to~$K_1$ without marked
vertices. Therefore, when we mark all vertices of $\^K_2$, we
obtain that any subgraph of $\Gamma$ isomorphic to $K_1$ has a marked
vertex. Thus, $\repV(K_1,\Gamma)\le\repV(K,\Gamma)+k_2$ and
$$
\repVSym(K,\Gamma)
\le
\repVSym(K_1,\Gamma)
\le
k_1\repV(K_1,\Gamma)
\le
k_1(\repV(K,\Gamma)+k_2)
\hbox{ (where the next to last inequality is Corollary 1)}
$$
as required.
If $\Gamma$ has no subgraphs isomorphic to $K_2$ without
marked vertices, then we obtain the inequality
$$
\repV(K_2,\Gamma)\le\repV(K,\Gamma)
$$
(which actually is an
equality) and
$$
\repVSym(K,\Gamma)
\le
\repVSym(K_2,\Gamma)
\le
k_2\repV(K_2,\Gamma)
\le
k_2\repV(K,\Gamma)
\le
k_1\repV(K,\Gamma)
<
k_1(\repV(K,\Gamma)+k_2)
$$
as required.
This completes the proof.

\Question 1.
Is any finite connected graph costly in the class of
connected graphs?

According to Theorem 1, all finite connected graphs without vertices
of degree one
are costly in the class of
connected graphs. Chains are also
costly in the class of
connected graphs; we can take polygons (cycles) as $\Gamma_m$ in this case.

The four-vertex graphs
\emph{tailed triangle}
and \emph{claw}
shown in Figure 2 on the left are also costly in the class of
connected graphs;
the corresponding graphs $\Gamma_m$ are shown in Figure 2
on the right. To be more precise, this figure shows an infinite
vertex-transitive graph where one-fourth of the vertices are marked in
such a way that each tailed triangle and each claw have a marked vertex;
this pattern on the plane is twice periodic, hence, we can obtain
arbitrarily large finite patterns on the torus with the same property
(i.e.  vertex-transitive graphs in which one quarter of the vertices are
marked, and each tailed triangle and each claw has a marked vertex).


\goodbreak
\bigskip
\centerline{\input 2.PIC}
\nobreak%
\centerline{Fig. \lowercase{2}}%
\goodbreak
\bigskip

Thus,
\disp{\sl
all
connected graphs
with at
most four vertices are costly in the class of
connected graphs.}%
Strictly speaking, Figure 2 proves this assertion only if
the word \emph{graph} means an undirected graph without loops and
multiple edges; but obvious modifications of this figure makes it
possible to prove the assertion in other cases.
For instance, if \emph{graph} means digraph without loops and multiple 
edges, than we should give some directions to the edges of the
tailed triangle and claw (Figure 2, left) and replace each edge by a pair 
of inverse edges in the graph shown in Figure 2 on the right.

Nevertheless, we conjecture that the answer to Question 1 is
negative, and a counter-example is probably the five-vertex graph
$D_5$ shown on Figure 3 (i.e., hypothetically
$\repVSym(D_5,\Gamma)<5\repV(D_5,\Gamma)$ for any connected graph
$\Gamma$ with sufficiently large representativeness $\repV(D_5,\Gamma)$).


\goodbreak
\bigskip
\centerline{\input 3.PIC}
\nobreak%
\centerline{Fig. \lowercase{3}}%
\goodbreak
\bigskip

A partial confirmation of this conjecture
is the following fact.

\Th 2.
The graph $D_5$ is not costly in
the
class vertex-transitive connected graphs.
More precisely,
if $\Gamma\supseteq D_5$ is a vertex-transitive undirected
connected graph
with more than five
vertices, and the representativeness $\repV(D_5,\Gamma)$
is finite, then
$\repVSym(D_5,\Gamma)<5\repV(D_5,\Gamma)$.

(A graph $\Gamma$ is called \emph{vertex-transitive} if its
automorphism group acts transitively on the set of vertices, i.e., for any
two vertices $u$ and $v$, there exists an automorphism $\phi$
of the graph such
that $\phi(u)=v$.)

\Proof
First note that we can (and shall) assume that $\Gamma$ has
no loops and multiple edges. Indeed, if we remove all loops, and
replace
every
bunch of
multiple edges with a single edge, then the conditions of the theorem
remain fulfilled, and the assertion for the obtained graph implies the
assertion for the initial graph.

Note also
that $\Gamma$ must be finite, because
by Corollary 1 the finiteness of $\repV(D_5,\Gamma)$ implies the finiteness
of $\repVSym(D_5,\Gamma)$, which is equal to the number
of vertices of the graph by virtue of transitivity.

Suppose that the degree of each vertex of $\Gamma$ is $k\ge3$
(if $k<3$, then we have nothing to prove).
Let us mark a finite set $X$ of vertices of $\Gamma$
(where $|X|=\repV(D_5,\Gamma)$)
in such a way that each subgraph isomorphic to $D_5$ contains
a marked vertex. This means that
the graph $\Gamma'$ obtained from $\Gamma$ by deleting all marked
vertices and edges incident to them contains no subgraphs isomorphic to
$D_5$. It is easy to classify such graphs.

\Lemma 1.
A finite connected
undirected graph without loops, multiple edges,
and subgraphs isomorphic to $D_5$ is either
\-
an $l$-gon (cycle) \(with $l\ge3$\),
\-
a chain \(with $l\ge0$ edges\),
\-
a star $K_{1,l}$ \(with $l\ge3$ edges\),
\-
or a connected four-vertex graph.

\noindent
Figure 4 shows (representatives of) the three infinite series;
Figure~5 shows three remaining four-vertex graph
(do not pay attention to black vertices and edges incident to them
for now).


\goodbreak
\bigskip
\centerline{\input 4.PIC}
\nobreak%
\centerline{Fig. \lowercase{4}}%
\goodbreak
\bigskip

\vfil\break

\goodbreak
\bigskip
\centerline{\input 5.PIC}
\nobreak%
\centerline{Fig. \lowercase{5}}%
\goodbreak
\bigskip

\Proof
If a connected finite graph contains no vertices of degree higher than
two, then this graph is either a polygon or a chain. If there is
a vertex $v$ of degree three or higher, then the
neighbouring vertices can be joined by edges only between themselves (and
$v$), because otherwise we obtain $D_5$ as a subgraph. Thus,
all vertices of the graph, except $v$, are neighbours of $v$ (by
virtue of connectedness).

If no neighbours of $v$ are joined by edges,
then the graph~ is a star.

If an edge joins two
neighbours $u$ and $w$ of $v$, then
at most one additional vertex
(except $v$, $u$, and $w$) can exist, because
otherwise we again obtain a subgraph isomorphic to $D_5$.
Thus, our graph has at most four vertices, that completes
the proof of the lemma.

\medskip

Proceeding with the proof of the theorem,
let us calculate the number $p$ of edges joining marked vertices with
non-marked ones. On the one hand,
$p\le k|X|$ (because each marked vertex has degree $k$),
and the equality is achieved only when no two marked
vertices are joined by an edge. On the other hand,
$p\ge(k-3)|Z|$, where $Z$ is the set of non-marked vertices, because
$$
{
\hbox{the number of edges joining a component of $\Gamma'$ with marked
vertices}
\over
\hbox{the number of vertices in this component}
}
\ge k-3
$$
(see Figs.~4~and~5, where marked vertices are black and
$k=4$), and the equality is achieved only on components which are
complete graphs with four vertices.

\noindent
Thus, $k|X|\ge p\ge(k-3)|Z|$, i.e. $|X|\ge(1-{3\over k})|Z|$.
This means that either
\item{1)}
$|X|>{1\over 4}|Z|$,
\item{2)}
$k=3$,
\item{3)}
or $k=4$ and  $|X|={1\over 4}|Z|$; as was mentioned, this
is possible only when all components of $\Gamma'$
are complete graphs on four vertices, and no two
marked vertices of $\Gamma$ are joined by an edge.

\enditem
Consider these cases.

\item{1)}
In this case, we have
$
\repVSym(D_5,\Gamma)\le
\hbox{(the number of vertices of $\Gamma$)}
=|X|+|Z|<|X|+4|X|=5|X|=5\repV(D_5,\Gamma)
$
as required.

\item{2)}
{%
In this case, a component of $\Gamma'$ cannot be a
complete graph on four vertices shown in Figure 5 on the left
(because the degree of each vertex of $\Gamma$ is three).
Neither can a component of $\Gamma'$ be a \emph{diamond},
shown on Figure 5 in the centre; indeed,
the neighbourhood of the vertex $u$ (in $\Gamma$) is a
chain (in this case), while the neighbourhood of the vertex $w$ is a
disconnected graph; this is impossible in a vertex-transitive graph.
Recall that the
\emph{neighbourhood} of a vertex $v$ is the graph consisting of
vertices neighbouring to $v$ and all edges between them. Note
that, in the case under consideration, the degree of each vertex is three
(not four, as on Figure 5).

Thus,
the argument that have led us to the inequality $k|X|\ge(k-3)|Z|$
are modified as follows:
$$
{
\hbox{the number of edges joining a component of $\Gamma'$ with marked
vertices}
\over
\hbox{the number of vertices in this component}
}
\ge 1
$$
(where the equality is achieved on components
shown in Figure~5 on the right and on cycles, Figure~4).
This implies that $3|X|\ge|Z|$ and we come to case 1).
}

\item{3)}
{%
The neighbourhood of a marked vertex is a
disjoint union of several cliques (consisting of
non-marked vertices).
By virtue of transitivity, the neighbourhood of a non-marked
vertex has the same form. Since the neighbourhood of a non-marked
vertex must contain a triangle consisting of non-marked vertices,
we come to the conclusion that either
\itemitem{a)}
the neighbourhood of each vertex is a complete graph on four
vertices
\itemitem{b)}
or no triangle contains a marked vertex.

In Case a) the graph $\Gamma$ must be the complete graph
with five vertices that completes the proof. Case b)
is impossible in a vertex-transitive graph
containing a cliques of order four.
}

\s 2.
Edge representativeness

We say that the \emph{edge representativeness} $\repE(K,\Gamma)$
of a graph $K$ in a graph $\Gamma$
is the minimal integer $n$ such that
$\Gamma$ contains a set $X$
of $n$ edges
satisfying the following condition:
$$
\hbox{\sl each subgraph of $\Gamma$ isomorphic to $K$ contains an edge
from $X$.}
\eqno{({**}*)}
$$
Let us define the
\emph{symmetric edge representativeness}
$\repESym(K,\Gamma)$ of a graph $K$ in a graph $\Gamma$ as the minimal
integer $n$ such that $\Gamma$ contains an $(\Aut\Gamma)$-invariant set
$X$ of $n$ edges satisfying $({**}*)$.

Clearly, $\repE(K,\Gamma)\le\repESym(K,\Gamma)$. Corollary 1 says that
$\repESym(K,\Gamma)\le
\repE(K,\Gamma)\cdot(\hbox{the number of edges of $K$})$.
We call a graph $K$
\emph{costly in the sense of symmetric edge representativeness}
or simply
\emph{edge-costly} if
$$
\hbox{
$\forall m\in\Z$ there exists a graph
$\Gamma_m$ such that
}
\repESym(K,\Gamma_m)=
\repE(K,\Gamma_m)\cdot(\hbox{the number of edges of $K$})\ge m.
\eqno{({**}{**})}
$$
Thus, a graph $K$ is edge-costly if the estimate from Corollary 1(2)
is sharp for $K$. We call an edge-costly
graph $K$ \emph{edge-costly in the class of graphs $\cal K$} if
the graphs $\Gamma_m$ in~$({**}{**})$ can chosen from the class $\cal K$.

\Proposition 1.
Any finite edge-transitive
connected
graph~$K$ is edge-costly in the class of connected graphs.

\Proof
If $K$ has no hanging edges (i.e. no vertices of degree
one), then we can do pretty much the same as in the vertex case.
Take as the graph~$\Gamma_m$
the disjoint union of $m$ copies of~$K$.
This shows that $K$ is edge-costly (in the class of all graphs yet).
Indeed, to destroy all subgraphs isomorphic to $K$,
it is
sufficient (and necessary) to remove
one edge from each copy of $K$. Thus,
$\repE(K,\Gamma_m)=m$.
The graph $\Gamma_m$ is edge-transitive, therefore,
$\repESym(K,\Gamma_m)=mk$, where $k$ is the number of edges of $K$.

To
make the graph $\Gamma_m$ connected,
we add to
$\Gamma_m$ chains of length~$N$ joining each
vertex of the $i$-th copy of~$K$ with the corresponding vertex of
the $(i+1)$-th
copy of~$K$, where $i\in\{1,\dots,m-1\}$, and the integer~$N$ (the
same for all chains) is larger than the number of vertices of~$K$.

Again, $\repE(K,\Gamma_m)=m$, because to destroy
the subgraphs isomorphic to $K$ we can remove one
edge 
from each copy of $K$ (as $K$ has no
hanging edges, and the number $N$ is large enough). Now,
$\repESym(K,\Gamma_m)=mk$, because
the group $\Aut\Gamma_m$ acts transitively on the edges of each copy
of~$K$.

If the edge-transitive connected graph $K$ has hanging edges, then all
edges are hanging and $K$ is a star.

If the edges are directed, then (by the edge-transitivity)
the star $K=K_{1,l}$ has one source and $l$ sinks or vice versa.
Assuming that there is one source, take the complete bipartite
graph $K_{m,l}=\Gamma_m$ in which all edge are directed from the first
part (consisting of $m$ vertices) to the second part (consisting of
$l$ vertices). Clearly, $\repE(K,\Gamma_m)=m$
and $\repESym(K,\Gamma_m)=ml$ (because the graph $\Gamma_m$
is edge-transitive).

It remains to consider the case, where $K=K_{1,l}$ is an
undirected star.

If $l$ equals one or two, then we can take the cycle of length $2m$
as
$\Gamma_m$.

If $l=3$, then we can take ``honeycombs" as
$\Gamma_m$.
``Honeycombs" (Fig. 6) is an
infinite graph, in which one third of the edges are marked
(vertical edge on Figure 6),
and each subgraph isomorphic to the claw $K=K_{1,3}$ has
a marked edge.
To make this graph finite, we note that the
``honeycombs" is a doubly periodic pattern on the plane; therefore,
we can
draw an arbitrarily large finite graph on the torus with the same
properties (i.e. an edge-transitive graph, in which one third of the edges
are marked, and each claw has a marked edge).

If $l>3$, then we can act
as shown in Figure 6 on the right: $\Gamma_m$ consists of $m$ copies of
the star $K$ joined by edges; $\repE(K,\Gamma_m)=m$ and
$\repESym(K,\Gamma_m)=ml$.
This completes the proof.

\vfil\break

\goodbreak
\bigskip
\centerline{\input 6.PIC}
\nobreak%
\centerline{Fig. \lowercase{6}}%
\goodbreak
\bigskip

The exact analogue of Theorem 1 is not true:
\disp{\sl
there exist disconnected edge-costly
graphs and even disconnected graphs edge-costly in the class of connected
graphs.
}%
Indeed, if take the disjoint union of
two edges as $K$, and the complete biparite graph $K_{2,m}$ as $\Gamma$,
then $\repE(K,\Gamma)=m$ (because, to destroy all
subgraphs of $\Gamma$ isomorphic to $K$, we can remove all edges
incident to a vertex of degree $m$, see Figure 7),
and $\repESym(K,\Gamma)=2m$ (because the graph is edge-transitive).

\goodbreak
\bigskip
\centerline{\input 7.PIC}
\nobreak%
\centerline{Fig. \lowercase{7}}%
\goodbreak
\bigskip

Are all graphs edge-costly?
No, but we have only trivial examples:
if the word ``graph'' means a directed graph
without multiple edges but
loops
are allowed, then the graph~$K$ shown in
Figure 8 is not edge-costly by a trivial reason: a common edge of two
subgraphs of this form in any graph $\Gamma$ must be a loop,
therefore $\repESym(K,\Gamma)=\repE(K,\Gamma)$
(because, if we want to destroy all subgraphs of this form,
it is more efficient to remove only loops).

\goodbreak
\bigskip
\centerline{\input 8.PIC}
\nobreak%
\centerline{Fig. \lowercase{8}}%
\goodbreak
\bigskip

We do not know, e.g., an answer to the following question.

\Question 2.
Let the word \emph{graph}
mean an undirected graph without loops and multiple edges.
Does there exist a finite non-edge-costly graph?
Can such a graph be connected?
Does there exist a finite graph which is not edge-costly in the class
of
connected graphs? Can such a graph be connected?

If the word \emph{graph}
means a directed graph without loops and multiple edges,
the similar questions are also open.

\s 3.
Proof of the main theorem

Put $m=\max\limits_{F\in\F}|F|$ and consider the following set
$
Y=\left\{y\in U\;\Bigm|\;
|Gy\cap X| \ge{1\over m}|Gy|\right\}
$
(in particular, $Y$ contains no points with infinite orbits).
Clearly, this set is $G$-invariant, and
$|Y|\le m|X|$ (because, for each orbit $Gu$, we have
$|Gu\cap Y|\le m|Gu\cap X|$).

It remains to show, that $Y$ is a system of representatives for $\F$.
Take a set $F\in\F$.
Each set $gF$ (where $g\in G$) belongs to $\F$ (because
the family $\F$ is $G$-invariant) and, hence, intersects $X$.
Therefore,
$$
G=\bigcup_{f\in F}\{g\in G\;|\;gf\in X\}.
$$
For 
each $f\in F$, the 
set $\{g\in G\;|\;gf\in X\}$ is
a union of finitely many (possibly zero) 
left cosets of the stabiliser $\St(f)$ of the point $f$; that is
$$
\{g\in G\;|\;gf\in X\}=
\bigcup_{x\in X}\{g\in G\;|\;gf=x\}=
\bigcup_{x\in X\cap Gf}g_x\cdot\St(f),
\qbox{where $g_x\in G$ are fixed such that $g_xf=x$}.
$$
Thus, we obtain a decomposition of the group~$G$ into a finite union
of left cosets of some subgroups.
B.~Neumann's theorem ([Neu54], Proposition~4.5) says that
\disp{\sl
if a group $G$ is covered by finitely many cosets of
subgroups: $G=g_1G_1\cup\dots\cup g_sG_s$,
then
\newline
$\displaystyle\sum {1\over|G:G_i|}\ge1$
\(where the inverse of an infinite cardinal
is zero\).
}
Therefore, (taking into account that
the index of a stabiliser equals the length of the corresponding orbit)
we obtain
$$
1\le\sum_{f\in F}{1\over|G:\St(f)|}\cdot|Gf\cap X|=
\sum_{f\in F}{|Gf\cap X|\over|Gf|}.
$$
Since the number of terms equals $|F|\le m$,
some term is at least $1/m$, i.e.,
${|Gf\cap X|/|Gf|}\ge1/m$, and, hence, $f\in Y$ (by the definition
of $Y$). This completes the proof.


\References

[BrNa04]
B. Bruno, F. Napolitani,
A note on nilpotent-by-\v Cernikov groups,
{Glasgow Math. J.}, 46 (2004), 211-215.

[ChD89]
A. Chermak, A. Delgado,
A measuring argument for finite group.
Proc. Amer. Math. Soc., 107 (1989), 907-914.

[dGT18a]
F. de Giovanni, M. Trombetti,
A note on large characteristic subgroups.
Communications in Algebra, 46:11 (2018), 4654-4662.

[dGT18b]
F. de Giovanni, M. Trombetti,
Large characteristic subgroups with modular subgroup lattice,
Archiv der Mathematik, 111:2 (2018), 123-128.

[dGT19a]
F. de Giovanni, M. Trombetti,
Large characteristic subgroups in which normality is a transitive relation,
Atti Accad. Naz. Lincei Cl. Sci. Fis. Mat. Natur. 30 (2019), 255-268.

[dGT19b]
F. de Giovanni, M. Trombetti,
Large characteristic subgroups and abstract group classes,
Quaestiones Mathematicae (to appear).

[Fr18]
E. Frolova,
Khukhro-Makarenko type theorems for algebras,
arXiv:1804.00268.

[Is08]
I. M. Isaacs,
Finite group theory,
GSM 92,
American Math. Soc.,
Providence RI, 2008.

[KaM82]
M. I. Kargapolov, Yu. I. Merzljakov,
Fundamentals of the theory of groups.
Graduate Texts in Mathematics, 62, Springer, 1979.

[KhM07a]
E. I. Khukhro, N. Yu. Makarenko,
Large characteristic subgroups satisfying multilinear commutator
identities,
J. London Math. Soc., 75:3 (2007), 635-646.

[KhM07b]
E. I. Khukhro, N. Yu. Makarenko,
Characteristic nilpotent subgroups of bounded co-rank and
\newline
automorphically-invariant ideals of bounded codimension in Lie algebras,
{Quart. J. Math.}, 58 (2007), 229-247.

[KhM08]
E. I. Khukhro, N. Yu. Makarenko,
Automorphically-invariant ideals satisfying multilinear identities,
and group-theoretic applications,
{J. Algebra}, 320:4 (2008), 1723--1740.

[KhKMM09]
E. I. Khukhro, Ant. A. Klyachko, N. Yu. Makarenko, and Yu. B. Melnikova
Automorphism invariance and identities.
Bull. London Math. Soc. (2009), 41(5): 804-816.
\arXiv 0812.1359

[KlMe09]
A. A. Klyachko, Yu. B. Mel'nikova,
A short proof of the Khukhro--Makarenko theorem on large characteristic
subgroups with laws,
Sbornik: Mathematics, 200:5 (2009), 661-664.
\arXiv 0805.2747

[KlMi15]
A. A. Klyachko, M. V. Milentyeva,
Large and symmetric: The Khukhro-Makarenko theorem on laws --- without laws,
J. Algebra, 424 (2015), 222-241.
\arXiv 1309.0571

[MSh12]
N. Yu. Makarenko, P. Shumyatsky,
Characteristic subgroups in locally finite groups,
J. Algebra, 352:1 (2012), 354-360.

[Neu54]
B. H. Neumann,
Groups covered by permutable subsets,
J. London Math. Soc., s1-29:2 (1954), 236-248.

[PSz02]
K. Podoski, B. Szegedy,
Bounds in groups with finite abelian coverings or with finite derived groups,
J. Group Theory, 5:4 (2002), 443-452.

[Vd00]
E. P. Vdovin,
Large normal nilpotent subgroups of finite groups,
Siberian Mathematical Journal, 41:2 (2000), 246-251.

\end